\documentclass[11pt]{amsart}
\usepackage{macros}

\title[Alternating super-elements of finite reflection groups]{Alternating super-polynomials and super-coinvariants \\ of finite reflection groups}
\author{Joshua P. Swanson}
\date{\today}

\begin{document}
  
\begin{abstract}
  Motivated by a recent conjecture of Zabrocki \cite{1902.08966}, Wallach
  \cite{1906.11787} described the alternants in the super-coinvariant
  algebra of the symmetric group in one set of commuting and one set of anti-commuting
  variables under the diagonal action. We give
  a type-independent generalization of Wallach's result to
  all real reflection groups $G$. As an intermediate step, we explicitly
  describe the alternating super-polynomials in $\bk[V] \otimes \Lambda(V)$
  for all complex reflection groups, providing an analogue of a classic result of
  Solomon \cite{MR0154929} which describes the invariant super-polynomials in
  $\bk[V] \otimes \Lambda(V^*)$. Using our construction, we explicitly
  describe the alternating harmonics and coinvariants for all real reflection groups.
\end{abstract}
\maketitle

\section{Introduction}

The classical coinvariant algebra of a complex reflection group $G \leq \GL(V)$
is the quotient $\bk[V]/\cI_+^G$ where $\bk[V] = \Sym(V^*)$ is
the $\bk$-algebra of polynomial functions on $V$, $\bk$ is a subfield of $\bC$,
and $\cI_+^G$ is the ideal generated by all homogeneous non-constant
$G$-invariants. Chevalley \cite{MR72877} showed that $\bk[V]/\cI_+^G$ as an
ungraded module carries the regular representation of $G$. The full graded
representation theory of coinvariant algebras and their generalizations is
extremely rich and has resulted in a vast body of work
(see e.g.~\cite{MR1168926,MR3783430,MR1256101,MR0354894,MR526968,MR1023791}).

Chevalley \cite{MR72877} and Shephard--Todd \cite{MR0059914}
further showed that $\bk[V]^G = \bk[f_1, \ldots, f_n]$ where
$f_1, \ldots, f_n$ are $n=\dim(V)$ homogeneous algebraically
independent $G$-invariants, and $\cI_+^G = \langle f_1, \ldots, f_n\rangle$.
Solomon generalized this result to the following explicit description of
the $G$-invariants of the Cartan algebra of differential forms on $V$,
$\bk[V] \otimes \Lambda(V^*)$. Here $d$ denotes the exterior derivative.

\begin{Theorem}[{Solomon \cite{MR0154929}}]\label{thm:inv_upstairs}
  For a complex reflection group $G$, the $G$-invariants of
  $\bk[V] \otimes \Lambda(V^*)$ have $\bk$-basis
  \begin{equation}
    \{ f_1^{\alpha_1} \cdots f_n^{\alpha_n} df_{i_1} \cdots df_{i_r}
        : \{i_1 < \cdots < i_r\} \subset [n], \alpha \in \bZ_{\geq 0}^n\}.
  \end{equation}
\end{Theorem}

\begin{Corollary}\label{cor:inv_upstairs_hilb}
  For a complex reflection group $G$ with degrees $d_1, \ldots, d_n$,
  \begin{equation}
    \Hilb((\bk[V] \otimes \Lambda(V^*))^G; q, t)
      = \prod_{i=1}^n \frac{1 + q^{d_i - 1} t}{1-q^{d_i}}.
  \end{equation}
\end{Corollary}

When $G = S_n$, the Cartan algebra $\bk[V] \otimes \Lambda(V^*)$
may be interpreted as the
ring of ``super-polynomials'' $\bk[\bm{x}_n, \bm{\theta}_n]$
in commuting variables $x_1, \ldots, x_n$ and anti-commuting
variables $\theta_1, \ldots, \theta_n$ where $\sigma(x_i) = x_{\sigma(i)},
\sigma(\theta_i) = \theta_{\sigma(i)}$. The $t=0$ specialization of a
recent conjecture of Zabrocki \cite{1902.08966}
concerns the bigraded $S_n$-module structure of the ``super-coinvariant
algebra'' $\bk[\bm{x}_n, \bm{\theta}_n]/\cJ_+$,
where $\cJ_+$ is the ideal generated by all homogeneous
$S_n$-invariant super-polynomials. Zabrocki's conjecture
provides an explicit module for the Delta conjecture
of Haglund--Remmel--Wilson \cite{MR3811519}, generalizing the relationship between
the diagonal coinvariants and the $n!$ theorem \cite{MR1214091}.
As a special case, Zabrocki's conjecture predicts that the bigraded Hilbert series
of the alternating component of
$\bk[\bm{x}_n, \bm{\theta}_n]/\cJ_+$ is
\begin{equation}\label{eq:Wallach}
  \Hilb((\bk[\bm{x}_n, \bm{\theta}_n]/\cJ_+)^{\det}; q, t)
    = \prod_{i=1}^{n-1} (t+q^i).
\end{equation}
Wallach \cite[Thm.~13]{1906.11787} has recently proven \eqref{eq:Wallach}.
Our primary objective is to give a type-independent generalization of \eqref{eq:Wallach}
valid for all real reflection groups.

Our first main result is an analogue of \Cref{thm:inv_upstairs}
for unitary reflection groups $G$ and the alternating component of
$\bk[V] \otimes \Lambda(V)$. Let $\Delta \in \bk[V]$ be the
generalized Vandermonde of $G$ (see \S\ref{ssec:vandermondes_jacobians})
and let $\odot$ be a certain action of $\bk[V] \otimes \Lambda(V^*)$ on
$\bk[V] \otimes \Lambda(V)$ involving partial differentiation and
multiplication operators (see \S\ref{sec:dots}).

\begin{Theorem}\label{thm:alt_upstairs}
  Let $G \leq U(n, \bk)$ be a unitary reflection group. Then the $\det$-isotypic component
  of $\bk[V] \otimes \Lambda(V)$ has $\bk$-basis
  \begin{equation}
    \{df_{i_1} \cdots df_{i_r} \odot (f_1^{\alpha_1} \cdots f_n^{\alpha_n} \Delta)
      : \{i_1 < \cdots < i_r\} \subset [n], \alpha \in \bZ_{\geq 0}^n\}.
  \end{equation}
\end{Theorem}

\begin{Corollary}\label{cor:alt_hilb}
  For a complex reflection group $G$ with degrees $d_1, \ldots, d_n$,
  \begin{equation}
    \Hilb((\bk[V] \otimes \Lambda(V))^{\det}; q, t)
      = \prod_{i=1}^n \frac{q^{d_i - 1} + t}{1 - q^{d_i}}.
  \end{equation}
\end{Corollary}

Solomon's description of $(\bk[V] \otimes \Lambda(V^*))^G$ motivates
the following natural generalization of Zabrocki's super-coinvariant algebra.

\begin{Definition}
  The \textit{super-coinvariant ideal} $\cJ_+^G$ is the ideal in $\bk[V] \otimes
  \Lambda(V^*)$ generated by bi-homogeneous non-constant $G$-invariant
  super-polynomials. The \textit{super-coinvariant algebra} of $G$ is
  $(\bk[V] \otimes \Lambda(V^*))/\cJ_+^G$. By Solomon's result,
  $\cJ_+^G = \langle f_1, \ldots, f_n, df_1, \ldots, df_n\rangle$.
\end{Definition}

\noindent The super-coinvariant algebra is a bigraded $G$-module. It
may be represented as a set of polynomials
$\cH_G \subset \bk[V] \otimes \Lambda(V^*)$, namely the
\textit{harmonics} of $\cJ_+^G$ relative to a certain
non-degenerate Hermitian form (see \S\ref{sec:herms_coinvs_harms}).
Precisely, the inclusion of $\cH_G$ into $(\bk[V] \otimes \Lambda(V^*))/\cJ_+^G$
is an isomorphism of bigraded $G$-modules. When $\bk \subset \bR$
and $G \leq O(n, \bk)$ is an orthogonal reflection group,
the representations $\Lambda(V)$ and $\Lambda(V^*)$ coincide,
and in this case we may explicitly describe the alternating component of
$\cH_G$ as follows.

\begin{Theorem}\label{thm:alt_harmonics}
  For a real reflection group $G \leq O(n, \bk)$ with fundamental invariants
  $f_1, \ldots, f_n$ of degrees $d_1 \leq \ldots \leq d_n$ with
  $f_1 = x_1^2 + \cdots + x_n^2$, the
  $\det$-isotypic component of $\cH_G$ has $\bk$-basis
  \begin{equation}
    \{df_{i_1} \cdots df_{i_r} \odot \Delta : \{i_1 < \cdots < i_r\} \subset [n]\}.
  \end{equation}
  Consequently, the projection of these elements yields a $\bk$-basis for
  $(\bk[V] \otimes \Lambda(V^*)/\cJ_+^G)^{\det}$.
\end{Theorem}

\begin{Corollary}\label{cor:alt_harmonics_hilb}
  For a real reflection group $G$ as in \Cref{thm:alt_harmonics},
  \begin{equation}
    \Hilb((\bk[V] \otimes \Lambda(V^*)/\cJ_+^G)^{\det}; q, t)
    = \prod_{i=1}^n (q^{d_i - 1} + t).
  \end{equation}
\end{Corollary}

Wallach's result \eqref{eq:Wallach} is the specialization of
\Cref{cor:alt_harmonics_hilb} where $G=S_n$ acts irreducibly on
$\bR^n/\langle 1, \ldots, 1\rangle$ via the standard representation
with degrees $2, \ldots, n$.

\begin{Remark}
  We may essentially always assume $G$ is unitary or orthogonal
  (see \S\ref{ssec:hermitian}). By choosing a suitable basis, we may
  further arrange for $x_1^2 + \cdots + x_n^2$ to be a fundamental
  invariant in the real case. The requirement $d_i \neq 1$ simply means
  $\{v \in V : \forall \sigma \in G, \sigma(v) = v\} = 0$,
  which is automatically satisfied for irreducible $G$. The assumptions
  in \Cref{thm:alt_upstairs} and \Cref{thm:alt_harmonics} are thus
  quite mild.
\end{Remark}

The rest of the paper is organized as follows. In \Cref{sec:preliminaries}
we give background on complex reflection groups. In \Cref{sec:molien}
we use Molien series to contrast different notions of super-alternants.
\Cref{sec:dots} describes two actions $\cdot$ and $\odot$ on super-polynomials
involving differential operators. \Cref{sec:herms_coinvs_harms}
generalizes standard terminology to super-coinvariants. Finally,
\Cref{sec:alt_upstairs} proves \Cref{thm:alt_upstairs} and
\Cref{sec:alt_harmonics} proves \Cref{thm:alt_harmonics}.

\section{Preliminaries}\label{sec:preliminaries}

Let $V$ be a finite-dimensional vector space over a field
$\bk \subset \bC$.

\subsection{Invariant Hermitian forms}\label{ssec:hermitian}

We first summarize some very well-known results concerning $G$-invariant
positive-definite Hermitian forms. We describe them in an unusual amount of
detail in order to make the choice of basis underlying
\Cref{thm:alt_harmonics} explicit.

\begin{Definition}
  A \textit{Hermitian form} on $V$ is a $\bZ$-bilinear map
  $\langle -, -\rangle \colon V \times V \to \bk$ such that
    \[ \langle \alpha u, \beta v\rangle = \overline{\alpha}\beta \langle u, v\rangle
        \qquad \text{and} \qquad
        \langle u, v\rangle = \overline{\langle v, u\rangle} \]
  for all $\alpha, \beta \in \bk$ and $u, v \in V$, where $\overline{\alpha}$
  is the complex conjugate of $\alpha$. Such a form is \textit{positive-definite}
  if additionally
  \begin{equation}
    \langle v, v\rangle > 0, \qquad \forall v \in V - \{0\}.
  \end{equation}
\end{Definition}

If $\bm{v}_i$ denotes the $i$th coordinate of $v \in V$ with respect to some
basis $B$, then $\sum_i \overline{\bm{u}_i} \bm{v}_i$ is the \textit{standard}
positive-definite Hermitian form on $V$ associated with $B$, so such forms
exist. If $G \leq \GL(V)$ is any finite subgroup and $\langle -, -\rangle'$ is
any positive-definite Hermitian form, it is easy to see that
  \[ \langle u, v\rangle \coloneqq \sum_{\sigma \in G} \langle \sigma(u), \sigma(v)\rangle' \]
is a $G$-invariant positive-definite Hermitian form, i.e.
  \[ \langle \sigma(u), \sigma(v)\rangle = \langle u, v\rangle \]
for all $\sigma \in G$. That is, $G$ consists of unitary transformations
with respect to $\langle -, -\rangle$.

Given a Hermitian form $\langle -, -\rangle$ on $V$ and an ordered basis $B$,
there is a unique matrix $A$ such that
  \[ \langle u, v\rangle = \bm{u}^\dagger A\bm{v} \]
where $\bm{v}$ is the column vector of $v$ with respect to $B$ and $\dagger$
denotes the conjugate-transpose. The matrix $A$ is \textit{Hermitian},
i.e.~$A^\dagger = A$. From the spectral theorem, $A$ is unitarily
diagonalizable, so $A = P^\dagger DP$ for $P^\dagger = P^{-1}$ and $D$
diagonal. Since $\langle -, -\rangle$ is positive-definite, the eigenvalues
of $A$ are positive reals. By extending scalars if necessary, we may
assume the square roots of the (positive, real) eigenvalues of $A$ are in $\bk$,
so we may write
$D = M^\dagger M$. Setting $Q \coloneqq MP$, we have $A = Q^\dagger Q$
where $Q$ is non-singular, so $\langle u, v\rangle = (Q\bm{u})^\dagger (Q\bm{v})$.
Consequently, we may replace $B$ with a new basis $B'$ for which $\langle u, v\rangle
= \bm{u}^\dagger \bm{v}$ is the standard Hermitian form. That is, we may
assume $G \subset U(n, \bk)$ consists of unitary matrices, so that $\sigma^{-1}
= \sigma^\dagger$ for all $\sigma \in G$. Furthermore, for the natural $G$-action on
the dual space $V^*$ with respect to the dual basis of $B'$,
$\sigma \in G \subset U(n, \bk)$ is
represented by $\overline{\sigma}$, the conjugate of $\sigma$, which remains unitary.

To summarize, we have the following. In practice, $G$ is often defined by generalized
permutation matrices, which are automatically unitary.

\begin{Lemma}\label{lem:Hermitian_invariant}
  For any finite subgroup $G$ of $\GL(V)$, possibly after extending scalars by square
  roots of positive reals, there exists
  a basis $B$ for which the actions of $G$ on $V$ and $V^*$ with respect to $B$
  and the dual basis of $B$ are unitary.
\end{Lemma}

\subsection{Complex reflection groups and super-polynomials}\label{ssec:super_polynomials}

\begin{Definition}
  A \textit{pseudoreflection} is an element $\sigma \in \GL(V)$ such that
    \[ \exists m \in \bZ_{\geq 1} \text{ s.t. } \sigma^m = 1 \qquad \text{and} \qquad 
        \codim\{v \in V : \sigma \cdot v = v\} = 1. \]
  A \textit{complex reflection group} is a finite subgroup $G \leq \GL(V)$
  generated by pseudoreflections.
\end{Definition}

\begin{Definition}
  Let $\bk[V] \coloneqq \Sym(V^*)$ be the ring of polynomial
  functions on $V$, namely the symmetric algebra on $V^*$ over
  $\bk$. If $V$ has basis $e_1, \ldots, e_n$ and $V^*$ has dual basis
  $x_1, \ldots, x_n$, we have
    \[ \bk[V] = \bk[x_1, \ldots, x_n] \eqqcolon \bk[\bm{x}_n]. \]
  The group $G \leq \GL(V)$ acts naturally on the dual $V^*$ via
    \[ (\sigma \cdot x)(v) \coloneqq x(\sigma^{-1}(v)) \]
  for all $\sigma \in G, x \in V^*, v \in V$.
  Similarly, $\bk[V]$ is naturally a graded $G$-module where
  $\sigma(fg) = \sigma(f)\sigma(g)$ for all
  $f, g \in \bk[V]$.
\end{Definition}

\begin{Definition}
  Chevalley \cite{MR72877} showed that the ring of polynomial $G$-invariants
  $\bk[V]^G$ is itself a polynomial ring generated by $n=\dim(V)$ homogeneous,
  algebraically independent elements $f_1, \ldots, f_n$ called
  \textit{fundamental invariants}, which are not unique. The multiset
  $d_1, \ldots, d_n$ of degrees of the fundamental invariants
  are the \textit{degrees} of $G$, which is unique.
\end{Definition}

\begin{Definition}
  Let $\Lambda(V^*)$ be the algebra of
  alternating multilinear functions on $V$ with values in $\bk$ under
  the wedge product, which can be realized as the exterior algebra of
  $V^*$ over $\bk$. $\Lambda(V^*)$ is naturally a graded $G$-module where
  $\sigma(f \wedge g) = \sigma(f) \wedge \sigma(g)$ for all $f, g \in \Lambda(V^*)$.
  To avoid confusion, we write
  $\theta_i$ instead of $x_i$ for the generators of $\Lambda(V^*)$,
  and we omit $\wedge$. Consequently, $\Lambda(V^*)$ has
  $\bk$-basis
    \[ \theta_I \coloneqq \theta_{i_1} \cdots \theta_{i_r} \]
  for $I = \{i_1 < \cdots < i_r\} \subset [n]$.
  By an abuse of notation, we write
    \[ \Lambda(V^*) = \bk[\theta_1, \ldots, \theta_n] \eqqcolon \bk[\bm{\theta}_n]. \]
\end{Definition}

\begin{Definition}
  Let $\Lambda(V)$ be the exterior algebra on $V$ over $\bk$,
  which is again a graded $G$-module. We write $\psi_i$ instead $e_i$
  for the generators of $\Lambda(V)$. As before, $\Lambda(V)$ has
  $\bk$-basis $\{\psi_I\}$ and we write
    \[ \Lambda(V) = \bk[\psi_1, \ldots, \psi_n] \eqqcolon \bk[\bm{\psi}_n]. \]
\end{Definition}

\begin{Definition}
  The \textit{super-polynomial rings} are $\bk[V] \otimes \Lambda(V^*)$
  and $\bk[V] \otimes \Lambda(V)$. We may write
  $\bk[V] \otimes \Lambda(V^*)$ as the associative $\bk$-algebra
    \[ \bk[x_1, \ldots, x_n, \theta_1, \ldots, \theta_n] \eqqcolon
        \bk[\bm{x}_n, \bm{\theta}_n] \]
  generated by indeterminates $x_1, \ldots, x_n$ and $\theta_1, \ldots, \theta_n$
  where
    \[ x_i x_j = x_j x_i \qquad x_i \theta_j = \theta_j x_i \qquad
        \theta_i \theta_j = -\theta_j \theta_i. \]
  Similarly we may realize $\bk[V] \otimes \Lambda(V)$ as
  $\bk[\bm{x}_n, \bm{\psi}_n]$. The super-polynomial rings are
  bigraded $G$-modules. We typically let $q$ track $\bm{x}$-degree
  and $t$ track $\bm{\theta}$ or $\bm{\psi}$-degree.
\end{Definition}

When $G$ consists of orthogonal matrices, we have $\Lambda(V) \cong \Lambda(V^*)$
via $\psi_i \mapsto \theta_i$ as $G$-modules, so in the real case there is
only one super-polynomial ring. We must be more careful in the complex case.

\begin{Definition}
  The \textit{differential} (or \textit{exterior derivative}) on $\bk[V] \otimes \Lambda(V^*)$
  is the $\bk$-linear map
  \begin{align*}
    d \colon \bk[\bm{x}_n, \bm{\theta}_n] &\to \bk[\bm{x}_n, \bm{\theta}_n] \\
    d(f(\bm{x}_n) \theta_I) &\coloneqq \sum_{i=1}^n \frac{\partial f(\bm{x}_n)}{\partial x_i}
      \theta_i \theta_I.
  \end{align*}
\end{Definition}

\noindent In particular, we have $dx_i = \theta_i$. The following observation is a routine
verification.

\begin{Lemma}[{\cite[p.~58]{MR0154929}}]\label{lem:d_G_action}
  For all $\sigma \in G$ and $f \in \bk[V] \otimes \Lambda(V^*)$,
  $\sigma \cdot d(f) = d(\sigma \cdot f)$.
\end{Lemma}

\begin{Remark}
  If $G$ consists of orthogonal matrices, then $x_1^2 + \cdots + x_n^2$ and
  $x_1 \theta_1 + \cdots + x_n \theta_n$ are both $G$-invariant elements of
  $\bk[V] \otimes \Lambda(V)$. If $G$ consists of unitary matrices, then
  $|x_1|^2 + \cdots + |x_n|^2$
  is $G$-invariant, which is not in general an element of either super-polynomial ring.
  However, $x_1 \psi_1 + \cdots + x_n \psi_n \in \bk[V] \otimes \Lambda(V)$ is
  $G$-invariant in general. For this reason, $\bk[V] \otimes \Lambda(V)$
  may be the more ``natural'' super-polynomial ring. Note, however, that
  Solomon's description of the $G$-invariants, \Cref{thm:inv_upstairs}, applies
  to $\bk[V] \otimes \Lambda(V^*)$ and not $\bk[V] \otimes \Lambda(V)$.
\end{Remark}

\subsection{Vandermondes and Jacobians}\label{ssec:vandermondes_jacobians}

Let $G$ be a complex reflection group with fixed fundamental
invariants $f_1, \ldots, f_n$.

\begin{Definition}
  The \textit{Vandermonde} of $G$ (relative to $f_1, \ldots, f_n$)
  is the element $\Delta \in \bk[V]$ defined by
  \begin{equation}
    df_1 \cdots df_n \eqqcolon \Delta \theta_{[n]}
      \in \bk[V] \otimes \Lambda(V^*).
  \end{equation}
  $\Delta$ may be thought of as the Jacobian determinant of
  $f_1, \ldots, f_n$ with respect to $x_1, \ldots, x_n$. Since
  $f_1, \ldots, f_n$ are algebraically independent, it follows that $\Delta \neq 0$.
\end{Definition}

\begin{Example}
  When $G=S_n$ and $f_i = \sum_j x_i^j$, we have
  $\Delta = c\prod_{i < j} (x_j - x_i)$ for $c = n!^n$.
\end{Example}

\begin{Remark}
  In fact, $\Delta$ has a very explicit description.
  Let $\ell_1, \ldots, \ell_r \in \bk[V]$ be linear forms which vanish on the $r$
  reflecting hyperplanes of $G$. Let $m_1, \ldots, m_r$ be the orders
  of the cyclic subgroups of $G$ fixing $\ker \ell_1, \ldots, \ker \ell_r$.
  Then (see \cite[p.~283]{MR0059914}, \cite{MR117285})
    \[ \Delta = c\prod_{i=1}^r \ell_i^{m_i - 1} \]
  for some $c \in \bk - \{0\}$.
\end{Remark}

\begin{Lemma}[{\cite[p.~59]{MR0154929}}]
  We have $\sigma(\Delta) = \det_V(\sigma) \Delta$ for all $\sigma \in G$.
\end{Lemma}

\section{Molien series and alternants}\label{sec:molien}

A classical theorem of Molien gives a succinct, beautiful, and
remarkably powerful description of the Hilbert series of the invariants
of the $G$-action on $\bk[V] = \bk[\bm{x}_n]$. An analogous result for
the $G$-action on $\Lambda(V^*) = \bk[\bm{\theta}_n]$ is less frequently
encountered, though the proof is no harder.
We require a bigraded generalization of these results
for relative invariants over fields other than $\bC$ for both types of super-polynomials.
Since it is difficult to find all of the relevant pieces for such a generalization in the
literature, we sketch a proof. We then use Molien's theorem, Solomon's
theorem, and a generalization of Solomon's theorem due to Orlik--Solomon
\cite{MR575083} to analyze several possible notions of alternants.

In this subsection only, $\bk$ denotes an arbitrary field of characteristic
zero, i.e.~not necessarily a subfield of $\bC$. For a $\bk G$-module $M$ and
$\sigma \in G$, let $\Tr_M(\sigma)$ denote the character of $M$ at $\sigma$,
and similarly let $\det_M(1 - \sigma q) \in \bk[q]$ denote the characteristic
polynomial of $\sigma$ acting on $M$. Note that $1/\det_M(1 - \sigma q) \in \bk[[q]]$
may be regarded as a formal power series and that $\bk[[q]] \supset \bQ[[q]]$.

\begin{Definition}
  If $M$ is irreducible, the \textit{Molien series} of a bigraded $\bk G$-module
  $\cS$ relative to $M$ is the formal power series
  \begin{equation}
    F_M(\cS; q, t) \coloneqq \sum_{i, j \geq 0} (\text{multiplicity of $M$ in the bidegree
      $(i, j)$ piece of $\cS$}) q^i t^j.
  \end{equation}
\end{Definition}

\begin{Theorem}[Molien]\label{thm:molien}
  Let $G$ be a finite subgroup of $\GL(V)$ where $V$ is a finite-dimensional
  vector space over a field $\bk$ of characteristic $0$. Suppose $M$ is
  an irreducible $\bk G$-module. Then
  \begin{align}
    F_M(\bk[V] \otimes \Lambda(V^*); q, t)
      &= \frac{1}{|G| \dim_{\bk} \bD} \sum_{\sigma \in G} \Tr_M(\sigma)
        \frac{\det_V(1 + \sigma t)}{\det_V(1 - \sigma q)}, \\
    F_M(\bk[V] \otimes \Lambda(V); q, t)
      &= \frac{1}{|G| \dim_{\bk} \bD} \sum_{\sigma \in G} \Tr_M(\sigma)
        \frac{\det_V(1 + \sigma^{-1} t)}{\det_V(1 - \sigma q)},
  \end{align}
  where $\bD \coloneqq \End_{\bk G}(M)$ is the division ring of
  $\bk G$-linear endomorphisms of $M$.
  
  \begin{proof}
    The following claim is well-known:
    \begin{align}\label{eq:molien.0}
      \frac{\dim_{\bD}(M)}{|G|} \sum_{\sigma \in G} \Tr_M(\sigma^{-1}) \sigma
    \end{align}
    is the unique $G$-invariant projection operator onto $M$-isotypic components. Indeed,
    by the Artin--Wedderburn theorem, $\bk G \cong \bigoplus_{i=1}^r \Mat_{n_i}(\bD_i)$
    where $M_1, \ldots, M_r$ are the inequivalent irreducible $\bk G$-modules,
    $\bD_i \coloneqq \End_{\bk G}(M_i)$ is a right $\bk$-module, and
    $n_i = \dim_{\bD_i}(M_i)$. Consequently,
    $\Tr_{\bk G}(\sigma) = \sum_{i=1}^r \dim_{\bD_i}(M_i) \Tr_{M_i}(\sigma)$. Let
    $e_1, \ldots, e_r \in \bk G$ be the primitive central orthogonal idempotents, so
    that for any $\bk G$-module $N$, $e_i$ is the unique $G$-invariant projection
    from $N$ to the $M_i$-isotypic component of $N$. Suppose $e_i = \sum_{\sigma \in G}
    a_\sigma \sigma$. We have for each $\sigma \in G$,
    \begin{align*}
      a_\sigma |G|
        = \Tr_{\bk G}(e_i \sigma^{-1})
        = \sum_{j=1}^r \dim_{\bD_j}(M_j) \Tr_{M_j}(e_i \sigma_j^{-1})
        = \dim_{\bD_i} \Tr_{M_i}(\sigma_i^{-1}),
    \end{align*}
    which proves the claim.
    
    Taking traces, it follows that the $\bk$-dimension
    of the $M$-isotypic component of $\bk[\bm{x}_n, \bm{\theta}_n]_{i, j}$ is
    \begin{equation}\label{eq:molien.1}
      \frac{\dim_{\bD}(M)}{|G|} \sum_{\sigma \in G} \Tr_M(\sigma^{-1})
      \Tr_{\bk[\bm{x}_n, \bm{\theta}_n]_{i,j}}(\sigma).
    \end{equation}
    Dividing \eqref{eq:molien.1} by $\dim_{\bk} M = \dim_{\bD} M \cdot \dim_{\bk} \bD$
    and using 
    \begin{equation}\label{eq:molien.2}
      \sum_{i, j \geq 0} \Tr_{\bk[\bm{x}_n, \bm{\theta}_n]_{i, j}}(\sigma) q^i t^j
        = \frac{\det_{V^*}(1+\sigma t)}{\det_{V^*}(1-\sigma q)}
        = \frac{\det_V(1+\sigma^{-1} t)}{\det_V(1-\sigma^{-1} q)}.
    \end{equation}
    yields the first Molien series, and the second is similar. \eqref{eq:molien.2}
    is straightforward to prove by diagonalizing $\sigma$.
  \end{proof}
\end{Theorem}

\begin{table}[ht]
  \centering
  \begin{tabular}{c|c|c|c}
    \ & $\det_V$ & $G$-invariants & $\det_V^{-1}$ \\
    \toprule
    $\bk[V] \otimes \Lambda(V)$
      & $\prod_{i=1}^n \frac{q^{e_i} + t}{1 - q^{d_i}}$
      & $\prod_{i=1}^n \frac{1 + q^{e_i^*} t}{1 - q^{d_i}}$
      & ? \\
    \midrule
    $\bk[V] \otimes \Lambda(V^*)$
      & $q^\delta \prod_{i=1}^n \frac{q^{e_i^*} + t}{1 - q^{d_i}}$
      & $\prod_{i=1}^n \frac{1+q^{e_i} t}{1 - q^{d_i}}$
      & $\prod_{i=1}^n \frac{q^{e_i^*} + t}{1 - q^{d_i}}$ \\
    \midrule
    $\bk[V^*] \otimes \Lambda(V)$
      & $\prod_{i=1}^n \frac{q^{e_i^*} + t}{1 - q^{d_i}}$
      & $\prod_{i=1}^n \frac{1 + q^{e_i} t}{1 - q^{d_i}}$
      & $q^\delta \prod_{i=1}^n \frac{q^{e_i^*} + t}{1 - q^{d_i}}$ \\
    \midrule
    $\bk[V^*] \otimes \Lambda(V^*)$
      & ?
      & $\prod_{i=1}^n \frac{1+q^{e_i^*} t}{1 - q^{d_i}}$
      & $\prod_{i=1}^n \frac{q^{e_i} + t}{1 - q^{d_i}}$
  \end{tabular}
  \vspace{0.5em}
  \caption{Product formulas for $\Hilb((\bk[V] \otimes \Lambda(V))^{\det}; q, t)$
  and related series.}
  \label{tab:products}
\end{table}

\begin{Remark}
  For general $G$, there are eight reasonable
  notions of ``alternant super-polynomials'': we may use $\bk[V]$ or $\bk[V^*]$;
  $\Lambda(V)$ or $\Lambda(V^*)$; and $\det$ or $\det^{-1}$.
  In the real case, $V \cong V^*$ and $\det = \det^{-1}$, so all eight
  notions coincide. In the complex case, they may be genuinely
  different. \Cref{thm:molien} allows us to relate many of the
  corresponding Hilbert series along with series for the $G$-invariants
  by simple transformations. For example,
  $\det_V(\sigma) \det_V(1 +\sigma t) = t^n \det_V(1+\sigma^{-1} t^{-1})$
  gives
    \[ \Hilb((\bk[V] \otimes \Lambda(V))^{\det_V}; q, t)
        = t^n \Hilb((\bk[V] \otimes \Lambda(V^*))^G; q, t^{-1}). \]
  \Cref{tab:products} lists product formulas for six of the eight
  types of alternants and all four types of $G$-invariants.
  In \Cref{tab:products}, we may go down and right one
  spot by applying $F(t) \mapsto t^n F(t^{-1})$ starting from the first or third rows;
   we may go down two spots and right one spot by applying
  $G(q) \mapsto (-q)^{-n} G(q^{-1})$; and
  we may reflect through the middle using $\sigma \mapsto \sigma^{-1}$,
  which preserves the Hilbert series. These operations result in three
  orbits.
  
  The orbit containing $(\bk[V] \otimes \Lambda(V^*))^G$ has four elements
  and yields product formulas arising from Solomon's result, \Cref{cor:inv_upstairs_hilb}.
  The orbit containing $(\bk[V] \otimes \Lambda(V))^G$ has six elements and
  yields product formulas arising from Orlik--Solomon's generalization of Solomon's
  result \cite[Thm.~3.1]{MR575083}. The remaining orbit of two elements is
  not covered by these results. In \Cref{tab:products}, $e_i \coloneqq d_i - 1$ are the
  \textit{exponents} of $G$, $e_1^*, \ldots, e_n^*$ are the
  \textit{coexponents} of $G$ defined by
    \[ \Hilb((\bk[V]/\cI_+^G \otimes V^*)^G; q) \eqqcolon
        q^{e_1^*} + \cdots + q^{e_n^*}, \]
  and $\delta \coloneqq e_1 + \cdots + e_n - e_1^* - \cdots - e_n^*$.
\end{Remark}

\begin{Remark}
  \cite[Thm.~2.5.3]{MR1249931} gives the relative Molien series
  for $\bk[V]$ (though $\dim_{\bk}(M)$ should instead be
  $\dim_{\bD}(M)$). As Benson notes, the result can be generalized
  as-is to non-modular fields using Brauer characters.
  When $\Char \bk=0$ and $\bk = \overline{\bk}$, the Schur orthogonality
  relations may be easily deduced from \eqref{eq:molien.0}.
\end{Remark}

\section{Two differential operator actions}\label{sec:dots}

It is well-known that the multivariate polynomial ring acts on itself
by polynomial differential operators. For the super-polynomial rings
$\bk[V] \otimes \Lambda(V^*) = \bk[\bm{x}_n, \bm{\theta}_n]$ and
$\bk[V] \otimes \Lambda(V) = \bk[\bm{x}_n, \bm{\psi}_n]$, the
anti-commuting variables may act either as a form of partial
differentiation, as in \cite{1906.11787} and \cite{1906.03315},
or they may act by multiplication.
Here we define actions of these super-polynomial rings on each other
and summarize their relationship. Many of these facts
(or the special case when $G=S_n$) appear in \cite{1906.11787} and
\cite[\S5]{1906.03315}. Similar actions appear in \cite{RSS19}.

\begin{Definition}
  We have three flavors of $\bk$-linear endomorphisms
  \begin{align*}
    \partial_i^x &\colon \bk[\bm{x}_n] \to
      \bk[\bm{x}_n] \\
    \partial_i^\theta, m_i^\theta &\colon \bk[\bm{\theta}_n] \to \bk[\bm{\theta}_n] \\
    \partial_i^\psi, m_i^\psi &\colon \bk[\bm{\psi}_n] \to \bk[\bm{\psi}_n]
  \end{align*}
  given by partial differentiation with respect to $x_i$ and partial differentiation or
  multiplication with respect to $\theta_i$ or $\psi_i$:
  \begin{align*}
    \partial_i f(\bm{x}_n) &\coloneqq \frac{\partial f(\bm{x}_n)}{\partial x_i} \\
    \partial_i^\theta(\theta_{i_1} \cdots \theta_{i_r}) &\coloneqq
      \begin{cases}
        (-1)^{m-1} \theta_{i_1} \cdots \widehat{\theta_{i_m}} \cdots \theta_{i_r}
          & \text{if $i_1 < \cdots < i_r$ and $i = i_m$} \\
        0 & \text{otherwise}
      \end{cases} \\
    m_i^\theta(\theta_I) &\coloneqq \theta_i \theta_I,
  \end{align*}
  and similarly with $\partial_i^\psi$ and $m_i^\psi$. We extend
  $\partial_i^x$ to yield $\bk[\bm{\theta}_n]$- or $\bk[\bm{\psi}_n]$-linear
  endomorphisms of the super-polynomial rings $\bk[\bm{x}, \bm{\theta}_n]$
  and $\bk[\bm{x}_n, \bm{\psi}_n]$. We similarly extend
  $\partial_i^\theta, m_i^\theta, \partial_i^\psi, m_i^\psi$ $\bk[\bm{x}_n]$-linearly
  to the super-polynomial rings.
\end{Definition}

\begin{Lemma}\label{lem:comm_acomm}
  We have the following commutation and anti-commutation relations.
  \begin{enumerate}[(i)]
    \item $[\partial_i^x, \partial_j^x] = 0$,
      $[\partial_i^x, \partial_j^\theta] = [\partial_i^x, m_j^\theta] = 0$,
      and $[\partial_i^x, \partial_j^\psi] = [\partial_i^x, m_j^\psi] = 0$.
    \item $\partial_i^\theta \partial_j^\theta = -\partial_j^\theta \partial_i^\theta, \quad
      m_i^\theta m_j^\theta = -m_j^\theta m_i^\theta, \quad
      \partial_i^\psi \partial_j^\psi = -\partial_j^\psi \partial_i^\psi, \quad
      m_i^\psi m_j^\psi = -m_j^\psi m_i^\psi$.
    \item $m_i^\theta \partial_j^\theta + \partial_j^\theta m_i^\theta = \delta_{i,j}, \quad
      m_i^\psi \partial_j^\psi + \partial_j^\psi m_i^\psi = \delta_{i,j}$,
  \end{enumerate}
  where $\delta_{i,j}$ is the identity if $i=j$ and $0$ otherwise.

  \begin{proof}
    In each case the identities involving $\psi$ are equivalent to the identities
    involving $\theta$, so we focus on the latter.
    \begin{enumerate}[(i)]
      \item The first equality is essentially classical and the second and third are immediate
        since $\partial_i^x$ operates on $x$-variables and $\partial_i^\theta, m_i^\theta$
        operate on $\theta$-variables,.
      \item The first equality is straightforward to verify on $\theta_I$ directly and
        extends $\bk[\bm{x}_n]$-linearly to all of $\bk[\bm{x}_n, \bm{\theta}_n]$.
        The second is immediate from $\theta_i\theta_j = -\theta_j\theta_i$.
      \item This is a consequence of the Leibniz rule, \Cref{lem:leibniz}, which we
        will prove shortly.
    \end{enumerate}
  \end{proof} 
\end{Lemma}

The following actions are fundamental to the rest of our arguments.

\begin{Definition}
  \ 
  \begin{enumerate}[label=(\alph*), ref=\theDefinition(\alph*)]
    \item\label{eq:cdot}
      We have an action of $\bk[V] \otimes \Lambda(V^*) = \bk[\bm{x}_n, \bm{\theta}_n]$
      on itself given by
      \begin{equation*}
        f(x_1, \ldots, x_n, \theta_1, \ldots, \theta_n) \cdot g \coloneqq
          \overline{f}(\partial_1^x, \ldots, \partial_n^x, \partial_1^\theta, \ldots,
          \partial_n^\theta)(g)
      \end{equation*}
      for all $f, g \in \bk[\bm{x}_n, \bm{\theta}_n]$, which is well-defined
      by \Cref{lem:comm_acomm}.
    \item\label{eq:odot} We have an action of
      $\bk[V] \otimes \Lambda(V^*) = \bk[\bm{x}_n, \bm{\theta}_n]$
      on $\bk[V] \otimes \Lambda(V) = \bk[\bm{x}_n, \bm{\psi}_n]$:
      \begin{equation*}
        f(x_1, \ldots, x_n, \theta_1, \ldots, \theta_n) \odot g \coloneqq
          \overline{f}(\partial_1^x, \ldots, \partial_n^x, m_1^\psi, \ldots, m_n^\psi)(g)
      \end{equation*}
      for all $f \in \bk[\bm{x}_n, \bm{\theta}_n], g \in \bk[\bm{x}_n, \bm{\psi}_n]$.
  \end{enumerate}
\end{Definition}

We emphasize the appearance of the coefficient-wise complex conjugate of $f$
in \Cref{eq:cdot} and \Cref{eq:odot}. This will be justified by an
equivariance property, \Cref{thm:unitary_actions}, which is our next goal.
We must first review the Leibniz rule for
$\partial_i^\theta$, which was stated in \cite{1906.03315},
though it was not used and the proof was left as an exercise to the reader. Since
the proof is somewhat intricate and we require \eqref{eq:leib.1} in
an essential way, we include a proof here.

\begin{Definition}
  For $I, J \subset [n]$, let
  \begin{equation}
    \inv(I, J) \coloneqq \#\{(i, j) \in I \times J : j < i\}.
  \end{equation}
  We see that if $I \cap J = \varnothing$,
    \[ \theta_I \theta_J = (-1)^{\inv(I, J)} \theta_{I \sqcup J}. \]
\end{Definition}

\begin{Lemma}[{\cite[(5.4)]{1906.03315}}]\label{lem:leibniz}
  For all $f, g \in \bk[\bm{x}_n, \bm{\theta}_n]$ where
  $f$ has $\theta$-degree $r$, we have the Leibniz rule
  \begin{equation}\label{eq:leib.1}
    \partial_i^\theta(fg) = \partial_i^\theta(f) g + (-1)^r f \partial_i^\theta(g).
  \end{equation}

  \begin{proof}
    We may suppose $f = \theta_I, g = \theta_J$
    where $I = \{i_1 < \cdots < i_r\}, J = \{j_1 < \cdots < j_s\}$. If 
    $I \cap J \neq \varnothing$, the left-hand side is $0$. If $I \cap J \supsetneq \{i\}$,
    both terms on the right-hand side are also $0$, so suppose $I \cap J = \{i\}$.
    Suppose $i = i_\ell = j_m$. The right-hand side is then
      \[ (-1)^{\ell-1} \theta_{I - \{i\}} \theta_J + (-1)^r \theta_I (-1)^{m-1} \theta_{J - \{i\}}.
      \]
    The powers on $-1$ differ by $(\ell - 1) - (r + m - 1)$.
    Now $\theta_{I - \{i\}} \theta_J$ and $\theta_I \theta_{J - \{i\}}$ differ
    in that $\theta_i$ has been commuted past $\theta_{i_{\ell+1}} \theta_{i_r}
    \theta_{j_1} \cdots \theta_{j_{m-1}}$, a total of $m-1 + r - \ell$ terms. It follows
    that the two terms are negatives, so they cancel. Thus, we may assume
    $I \cap J = \varnothing$. If $i \not\in I \cup J$,
    then each term is $0$, so we may suppose $i \in I \cup J$. The left-hand side is then
      \[ \partial_i^\theta (-1)^{\inv(I, J)} \theta_{I \cup J}
          = (-1)^{\#I \cup J < i + \#J < I} \theta_{I \cup J - \{i\}}, \]
    where $\#I \cup J < i$ is shorthand for the number of elements of $I \cup J$
    smaller than $i$, and similarly with $\#J < I$.
    
    If $i \in I$, then $i \not\in J$, and the right-hand side becomes
    \begin{align*}
      (-1)^{\#I < i} \theta_{I - \{i\}} \theta_J + 0
        &= (-1)^{\#I < i} (-1)^{\inv(I - \{i\}, J)} \theta_{I \cup J - \{i\}} \\
        &= (-1)^{\#I < i + \#J < I - \{i\}} \theta_{I \cup J - \{i\}}.
    \end{align*}
    Now
      \[ \#I \cup J < i + \#J < I \equiv_2 \#I < i + \#J < I - \{i\} \]
    is equivalent to
      \[ \#J < i + \#J < i \equiv_2 0, \]
    which is true.
    
    On the other hand, if $i \in J$, then $i \not\in I$, and the
    right-hand side becomes
    \begin{align*}
      0 + (-1)^{\#I} (-1)^{\#J < i} \theta_I \theta_{J - \{i\}}
        &= (-1)^{\#I} (-1)^{\#J < i} (-1)^{\inv(I, J-\{i\})} \theta_{I \cup J - \{i\}} \\
        &= (-1)^{\#I + \#J < i + \#J - \{i\} < I} \theta_{I \cup J - \{i\}}.
    \end{align*}
    We see directly that
      \[ \# I \cup J < i + \#J < I = \#I + \#J < i + \#J - \{i\} < I \]
    is equivalent to
      \[ \#I < i + \# i < I = \#I, \]
    which is true since $i \not\in I$. This completes the proof of \eqref{eq:leib.1}.
  \end{proof}
\end{Lemma}

\begin{Lemma}\label{lem:unitary_conj}
  Suppose $\sigma \in U(n, \bk)$ is unitary and $x_1, \ldots, x_n$
  are the coordinate functions on $\bk^n$. Then for all
  $g \in \bk[\bm{x}_n, \bm{\theta}_n]$ and $h \in \bk[\bm{x}_n, \bm{\psi}_n]$,
  \begin{align}
    \label{eq:conj.1} (\sigma \circ \partial_i^x \circ \sigma^{-1})(g) = \sigma(x_i) \cdot g
      \qquad&\text{and}\qquad
      (\sigma \circ \partial_i^x \circ \sigma^{-1})(h) = \sigma(x_i) \odot h \\
    \label{eq:conj.2} (\sigma \circ \partial_i^\theta \circ \sigma^{-1})(g) &= \sigma(\theta_i) \cdot g \\
    \label{eq:conj.3} (\sigma \circ m_i^\psi \circ \sigma^{-1})(h) &= \sigma(\theta_i) \odot h.
  \end{align}
  
  \begin{proof}
    For \eqref{eq:conj.1}, by $\bk[\bm{\theta}_n]$- or $\bk[\bm{\psi}_n]$-linearity,
    it suffices to consider the case when
    $g = x_1^{\alpha_1} \cdots x_n^{\alpha_n}$. Let
    $\sigma(u_i) \coloneqq x_i$, so $u_1, \ldots, u_n$ forms another basis of $V^*$.
    Suppose $x_i = \sum_j c_{ij} u_j$ and $u_i = \sum_j d_{ij} x_j$, so that
    $\frac{\partial x_i}{\partial u_j} = c_{ij}$ and $\frac{\partial u_i}{\partial x_j} = d_{ij}$.
    Furthermore, $\sigma(x_i) = \sigma(\sum_j c_{ij} u_j) = \sum_j c_{ij} x_j$, so
    $[c_{ij}]$ is the matrix of $\sigma$, and $\sigma^{-1}(x_i) = u_i = \sum_j d_{ij} x_j$,
    so $[d_{ij}]$ is the matrix of $\sigma^{-1}$. Since $\sigma$ is assumed unitary,
    we have $[c_{ij}] = [d_{ij}]^\dagger$, so $c_{ij} = \overline{d_{ji}}$.
    Using the multivariate chain rule, we now compute
    \begin{align*}
      (\sigma \circ \partial_i^x \circ \sigma^{-1})(x_1^{\alpha_1} \cdots x_n^{\alpha_n})
        &= \sigma\left(\frac{\partial}{\partial x_i}u_1^{\alpha_1} \cdots u_n^{\alpha_n}\right) 
        = \sigma\left(\sum_{j=1}^n \frac{\partial}{\partial u_j} u_1^{\alpha_1} \cdots u_n^{\alpha_n} \frac{\partial u_j}{\partial x_i}\right) \\
        &= \sigma\left(\sum_{j=1}^n \alpha_j u_1^{\alpha_1} \cdots u_j^{\alpha_j - 1} \cdots u_n^{\alpha_n} d_{ji}\right) \\
        &= \sum_{j=1}^n \alpha_j x_1^{\alpha_1} \cdots x_j^{\alpha_j - 1} \cdots x_n^{\alpha_n} \overline{c_{ij}} \\
        &= \sum_{j=1}^n \frac{\partial}{\partial x_j} x_1^{\alpha_1} \cdots x_n^{\alpha_n} \overline{c_{ij}}
        = \left(\sum_{j=1}^n c_{ij} x_j \right) \cdot x_1^{\alpha_1} \cdots x_n^{\alpha_n} \\
        &= \sigma\left(\sum_{j=1}^n c_{ij} u_j\right) \cdot x_1^{\alpha_1} \cdots x_n^{\alpha_n}
        = \sigma(x_i) \cdot x_1^{\alpha_1} \cdots x_n^{\alpha_n}.
    \end{align*}
    This proves \eqref{eq:conj.1}.
    
    For \eqref{eq:conj.2}, we begin with an analogue of the multivariate
    chain rule in this context. Let $\phi_i \coloneqq du_i$, so
    $\phi_1, \ldots, \phi_n$ is a basis for $\Lambda^1(V^*)$. We claim that
    \begin{equation}\label{eq:anticomm_chain}
      \frac{\partial}{\partial \theta_i} \phi_{k_1} \cdots \phi_{k_r}
        = \sum_{j=1}^r \frac{\partial}{\partial \phi_{k_j}} \phi_{k_1}
           \cdots \phi_{k_r} \frac{\partial \phi_{k_j}}{\partial \theta_i}
    \end{equation}
    for all $k_1 < \cdots < k_r$. When $r=0$, the result is clear. For $r>0$,
    by induction and \eqref{eq:leib.1}, we have
    \begin{align*}
      \frac{\partial}{\partial \theta_i} \phi_{k_1} \cdots \phi_{k_r}
        &= \frac{\partial \phi_{k_1}}{\partial \theta_i} \phi_{k_2} \cdots \phi_{k_r}
          - \phi_{k_1} \frac{\partial}{\partial \theta_i} \phi_{k_2} \cdots \phi_{k_r} \\
        &= \phi_{k_2} \cdots \phi_{k_r} \frac{\partial \phi_{k_1}}{\partial \theta_i} 
          - \phi_{k_1} \sum_{j=2}^r \frac{\partial}{\partial \phi_{k_j}} \phi_{k_2}
            \cdots \phi_{k_r} \frac{\partial \phi_{k_j}}{\partial \theta_i} \\
        &= \widehat{\phi_{k_1}} \phi_{k_2} \cdots \phi_{k_r}
          \frac{\partial \phi_{k_1}}{\partial \theta_i} 
          - \sum_{j=2}^r (-1)^{j-2} \phi_{k_1} \phi_{k_2} \cdots \widehat{\phi_{k_j}}
            \cdots \phi_{k_r} \frac{\partial \phi_{k_j}}{\partial \theta_i} \\
        &= \sum_{j=1}^r \frac{\partial}{\partial \phi_{k_j}} \phi_{k_1} \cdots \phi_{k_r}
          \frac{\partial \phi_{k_j}}{\partial \theta_i},
    \end{align*}
    proving \eqref{eq:anticomm_chain}. Now \eqref{eq:conj.2} follows from
    virtually the same calculation as \eqref{eq:conj.1} using \eqref{eq:anticomm_chain};
    the details are omitted.
    
    As for \eqref{eq:conj.3}, we have
    \begin{align*}
      (\sigma \circ m_i^\psi \circ \sigma^{-1})(g)
        &= \sigma(\psi_i \sigma^{-1}(g))
          = \sigma(\psi_i) g
          = \sum_j f_{ij} \psi_j g
          = \left(\sum_j \overline{f_{ij}} \theta_j\right) \odot g
    \end{align*}
    where $[f_{ij}]$ is the matrix of $\sigma$ with respect to $\psi_1, \ldots, \psi_n$.
    Since $\psi_1, \ldots, \psi_n$ is the dual basis of $\theta_1, \ldots, \theta_n$,
    the matrix of $\sigma$ acting on $V$ is the inverse-transpose of the matrix
    of $\sigma$ acting on $V^*$, namely
    $([f_{ij}]^{-1})^T = ([f_{ij}]^{\dagger})^T = [\overline{f_{ij}}]$. Thus
    $\sum_j \overline{f_{ij}} \theta_j = \sigma(\theta_i)$, completing the result.
  \end{proof}
\end{Lemma}

\begin{Theorem}\label{thm:unitary_actions}
  If $\sigma \in U(n, \bk)$ is unitary, then for all
  $f, g \in \bk[\bm{x}_n, \bm{\theta}_n]$ and
  $h \in \bk[\bm{x}_n, \bm{\psi}_n]$,
  \begin{align}
    \label{eq:dot_odot.1} \sigma(f \cdot g) &= \sigma(f) \cdot \sigma(g) \\
    \label{eq:dot_odot.2} \sigma(f \odot h) &= \sigma(f) \odot \sigma(h).
  \end{align}
  
  \begin{proof}
    For \eqref{eq:dot_odot.1}, by \eqref{eq:conj.1} and \eqref{eq:conj.2},
    \begin{align*}
      (\sigma \circ f(\partial_1^x, &\ldots, \partial_n^x, \partial_1^\theta,
      \ldots, \partial_n^\theta) \circ \sigma^{-1})(g) \\
        &= f(\sigma \circ \partial_1^x \circ \sigma^{-1}, \ldots,
          \sigma \circ \partial_n^x \circ \sigma^{-1}, \sigma \circ
          \partial_1^\theta \circ \sigma^{-1}, \ldots, \sigma \circ
          \partial_n^\theta \circ \sigma^{-1})(g) \\
        &= f(\sigma(x_1), \ldots, \sigma(x_n), \sigma(\theta_1), \ldots, \sigma(\theta_n))
          \cdot g \\
        &= \sigma(f) \cdot g.
    \end{align*}
    Replacing $g$ with $\sigma(g)$ gives the result. Similarly \eqref{eq:dot_odot.2}
    follows from \eqref{eq:conj.1} and \eqref{eq:conj.3}.
  \end{proof}
\end{Theorem}

\subsection{Hodge duality}

The two actions $\cdot$ and $\odot$ are related by the following
operation. We will not directly use the results of this subsection but
include it for completeness.

\begin{Definition}
  The \textit{Hodge dual} is the $\bk[\bm{x}_n]$-linear
  endomorphism
  \begin{align*}
    \hstar \colon \bk[\bm{x}_n, \bm{\theta}_n] &\to \bk[\bm{x}_n, \bm{\psi}_n] \\
    \hstar\theta_I &\coloneqq (-1)^{\deg(I)} \psi_J
  \end{align*}
  where $I \sqcup J = [n]$ and
    \[ \deg(I) \coloneqq \sum_{i \in I} (i-1). \]
\end{Definition}

\begin{Lemma}
  For $f, g \in \bk[\bm{x}_n, \bm{\theta}_n]$,
  we have
  \begin{equation}\label{eq:hodge_dots}
    f \odot \hstar g = \hstar(f \cdot g)
  \end{equation}
  
  \begin{proof}
    By $\bk[\bm{x}_n]$-sesquilinearity, we may suppose
    $f = \theta_I$ and $g = \theta_{[n]-J}$.
    The left-hand side becomes
      \[ \theta_I \odot (-1)^{\deg([n]-J)} \psi_{J}
          = (-1)^{\deg([n]-J)} \psi_I \psi_J
          = (-1)^{\deg([n]-J) + \inv(I, J)} \psi_{I \sqcup J} \]
    or $0$ if $I \cap J \neq \varnothing$. The right-hand side becomes
      \[ \hstar(\theta_I \cdot \theta_{[n]-J})
          = \hstar((-1)^{\inv(I, [n]-J)} \theta_{[n]-J-I})
          = (-1)^{\inv(I, [n]-J) + \deg([n]-J-I)} \psi_{I \cup J} \]
    or $0$ if $I \not\subset [n] - J$. Since $I \cap J = \varnothing$ if and only if
    $I \subset [n]-J$, we may suppose $I \cap J = \varnothing$.
    The result is hence equivalent to
    \begin{equation}\label{eq:hodge_dots.1}
      \deg([n] - J) + \inv(I, J) \equiv_2 \inv(I, [n]-J) + \deg([n]-I-J).
    \end{equation}
    Since $I \cap J = \varnothing$, we have $[n]-J = ([n]-J-I) \sqcup I$
    and $\deg([n]-J) = \deg([n]-J-I) + \deg(I)$. Thus \eqref{eq:hodge_dots.1}
    becomes
    \begin{equation}\label{eq:hodge_dots.2}
      \deg(I) + \inv(I, J) \equiv_2 \inv(I, [n]-J).
    \end{equation}
    Indeed, for all $I, J \subset [n]$, we have
    \begin{equation}\label{eq:hodge_dots.3}
      \inv(I, J) + \inv(I, [n]-J) = \deg(I)
    \end{equation}
    since
      \[ \#J < I + \#([n]-J) < I = \#[n] < I = \deg(I). \]
  \end{proof}
\end{Lemma}

\begin{Corollary}
  Let $f, g \in \bk[\bm{x}_n, \bm{\theta}_n]$ have the same bi-degree. Then
  \begin{equation}
    f \odot \hstar g = \langle f, g\rangle \psi_{[n]}
  \end{equation}
  where $\langle -, -\rangle$ is the non-degenerate Hermitian form
  defined in the next section.
\end{Corollary}

\section{Hermitian forms, coinvariants, and harmonics}\label{sec:herms_coinvs_harms}

We next define a non-degenerate Hermitian form on the super-polynomial
ring $\bk[V] \otimes \Lambda(V^*)$. We then summarize the connection
between the harmonic super-polynomials and super-coinvariants.

\begin{Definition}
  We have a $\bZ$-bilinear form on $\bk[V] \otimes \Lambda(V^*)
  = \bk[\bm{x}_n, \bm{\theta}_n]$ given by
  \begin{equation}
    \langle f, g\rangle \coloneqq \text{constant coefficient of $f \cdot g$}.
  \end{equation}
\end{Definition}

\begin{Lemma}\label{lem:herm_nondeg}
  The form $\langle -, -\rangle$ on $\bk[\bm{x}_n, \bm{\theta}_n]$
  is Hermitian, non-degenerate, and $G$-invariant. Moreover, for
  $\alpha, \beta \in \bZ_{\geq 0}^n$ and $I, J \subset [n]$,
  \begin{equation}\label{eq:orth}
    \langle\bm{x}^\alpha \theta_I, \bm{x}^\beta \theta_J\rangle
      = \begin{cases}
           (-1)^{\binom{|I|}{2}} \alpha! & \text{if $\alpha = \beta$, $I = J$} \\
            0 & \text{otherwise},
         \end{cases}
  \end{equation}
  where $\bm{x}^\alpha \coloneqq x_1^{\alpha_1} \cdots x_n^{\alpha_n}$,
  $\alpha! \coloneqq \alpha_1! \cdots \alpha_n!$. Consequently, $\langle -, -\rangle$
  is positive-definite or negative-definite when restricted to $\theta$-degree $r$
  depending on the sign of $(-1)^{\binom{r}{2}}$.
  
  \begin{proof}
    It is clear that $\langle f, g\rangle$ is conjugate-linear in the first argument
    and linear in the second argument. $G$-invariance follows from the first
    half of \Cref{thm:unitary_actions}. Non-degeneracy and the symmetry
    $\langle f, g\rangle = \overline{\langle g, f\rangle}$ both follow from
    \eqref{eq:orth}. As for \eqref{eq:orth}, we may assume that
    $\bm{x}^\alpha \theta_I$ and $\bm{x}^\beta \theta_J$ have
    the same bi-degree, in which case we see that
    \begin{align*}
      (x^\alpha \theta_I) \cdot (x^\beta \theta_J)
        &= (x^\alpha \cdot x^\beta)(\theta_I \cdot \theta_J) \\
        &= \alpha! \delta_{\alpha=\beta} (\theta_I \cdot \theta_I) \delta_{I=J}.
    \end{align*}
    It is straightforward to check that $\theta_I \cdot \theta_I
    = (-1)^{(|I|-1) + (|I|-2) + \cdots + 0} = (-1)^{\binom{|I|}{2}}$.
  \end{proof}
\end{Lemma}

\begin{Definition}
  Given a subspace $W$ in $\bk[\bm{x}_n, \bm{\theta}_n]$, the
  \textit{orthogonal complement} of $W$ is
  \begin{equation}
    W^\perp \coloneqq \{g \in \bk[\bm{x}_n, \bm{\theta}_n] : \langle f, g\rangle = 0
      \text{ for all }f \in W\}.
  \end{equation}
\end{Definition}

\begin{Lemma}
  We have
  \begin{enumerate}[(i)]
    \item $W^\perp
      = \{f \in \bk[\bm{x}_n, \bm{\theta}_n] : \langle f, g\rangle = 0
      \text{ for all }g \in W\}$
    \item $(W^\perp)^\perp = W$
    \item $W \cap W^\perp = 0$ and $W \oplus W^\perp = \bk[\bm{x}_n, \bm{\theta}_n]$
  \end{enumerate}
  
  \begin{proof}
    These are all standard consequences of \Cref{lem:herm_nondeg}.
  \end{proof}
\end{Lemma}

\begin{Definition}
  Let $G \leq U(n, \bk)$. The \textit{coinvariant ideal} of $G$ is
  the ideal
    \[ \cJ_+^G \coloneqq (\text{homogeneous non-constant $G$-invariants})
        \subset \bk[V] \otimes \Lambda(V^*). \]
  The \textit{coinvariant algebra} of $G$ is the quotient
    \[ \bk[V] \otimes \Lambda(V^*)/\cJ_+^G, \]
  which is a bigraded $G$-module. The \textit{$G$-harmonics} are
    \[ \cH_G \coloneqq (\cJ_+^G)^\perp \subset \bk[V] \otimes \Lambda(V^*). \]
\end{Definition}

The classical harmonic polynomials are those for which
$\sum_i \frac{\partial^2}{\partial x_i^2} f = 0$.
The $G$-harmonics are so named because of the following well-known result.

\begin{Lemma}\label{lem:harmonics_annihilated}
  If $\cJ \subset \bk[\bm{x}_n, \bm{\theta}_n]$ is an ideal generated by
  $j_1, \ldots, j_r$, then
  \begin{align*}
    \cJ^\perp
      &= \{g \in \bk[\bm{x}_n, \bm{\theta}_n] : j \cdot g = 0 \text{ for all }j \in \cJ\} \\
      &= \{g \in \bk[\bm{x}_n, \bm{\theta}_n] :
        j_1 \cdot g = \cdots = j_r \cdot g = 0\}.
  \end{align*}
  In particular, when $G$ is a complex reflection group with
  fundamental invariants $f_1, \ldots, f_n$,
    \[ \cH_G = \{g \in \bk[\bm{x}_n, \bm{\theta}_n] :
        f_1 \cdot g = \cdots = f_n \cdot g = df_1 \cdot g = \cdots = df_n \cdot g = 0\}. \]
  
  \begin{proof}
    If $j \cdot g= 0$ for all $j \in \cJ$, then trivially $\langle j, g\rangle = 0$.
    The converse also holds since $\cJ$ is closed under multiplication by
    $x_i$ and $\theta_i$. Similarly, $\cJ$ annihilates $g$ if and only if
    the generators annihilate $g$. The explicit description of $\cH_G$
    follows from Solomon's \Cref{thm:inv_upstairs}.
  \end{proof}
\end{Lemma}

\begin{Proposition}\label{prop:projection_iso}
  The projection of $\cH_G$ to $\bk[V] \otimes \Lambda(V^*)/\cJ_+^G$ is
  an isomorphism of bigraded $G$-modules.
  
  \begin{proof}
    $\cH_G$ is closed under the $G$-action since $\langle -, -\rangle$ and
    $\cJ_+^G$ are $G$-invariant. The projection is trivially $G$-equivariant
    and bidegree-preserving.
    It is a bijection since $\bk[\bm{x}_n, \bm{\theta}_n]/\cJ_+^G
    = ((\cJ_+^G)^\perp \oplus \cJ_+^G)/\cJ_+^G \equiv (\cJ_+^G)^\perp = \cH_G$.
  \end{proof}
\end{Proposition}

\section{Alternants in $\bk[V] \otimes \Lambda(V)$}\label{sec:alt_upstairs}

We may now prove the claimed classification of the alternants in
$\bk[V] \otimes \Lambda(V)$ for $G \leq U(n, \bk)$ from the introduction.

\begin{proof}[Proof of \Cref{thm:alt_upstairs}]
  Recall that the claimed $\bk$-basis for $(\bk[V] \otimes \Lambda(V))^{\det}$ is
    \[ \{df_I \odot f^\alpha \Delta
        : I \subset [n], \alpha \in \bZ_{\geq 0}^n\}, \]
  where $df_I \coloneqq df_{i_1} \cdots df_{i_r}$ for
  $I = \{i_1 < \cdots < i_r\}$ and $f^\alpha \coloneqq f_1^{\alpha_1} \cdots
  f_n^{\alpha_n}$. We first show that $\{df_I \odot f^\alpha \Delta\}$ carries the
  $\det$-representation and is linearly independent.
  
  As is well-known \cite{MR0154929}, $\bk[V]^{\det} = \bk[V]^G \Delta$,
  so $\{f^\alpha \Delta : \alpha \in \bZ_{\geq 0}^n\}$
  is a $\bk$-basis for $\bk[V]^{\det}$. Since $f_i \in \bk[V]^G$,
  we have $df_I \in (\bk[V] \otimes \Lambda(V^*))^G$ by
  \Cref{lem:d_G_action} and 
  $df_I \odot f^\alpha \Delta \in (\bk[V] \otimes \Lambda(V))^{\det}$
  by the second part of \Cref{thm:unitary_actions}.
  
  Now let $\{\tilde{g}_\alpha\}$ be
  an orthogonal basis for $\bk[V]^{\det}$ with respect to the
  non-degenerate Hermitian form $\langle -, -\rangle$. Since
  $\bk[V]^{\det} = \bk[V]^G \Delta$, we have $\tilde{g}_\alpha
  = g_\alpha \Delta$ for some $g_\alpha \in \bk[V]^G$. As far as
  linear independence is concerned, we may replace $f^\alpha \Delta$ with
  $g_\alpha \Delta$. Consequently, suppose
  \begin{equation}\label{eq:alt_upstairs.1}
    0 = \sum_{J, \beta} c_{J, \beta} df_J \odot g_\beta \Delta
  \end{equation}
  for some $c_{J, \beta} \in \bk$. By homogeneity in the $\bm{\theta}$
  variables, we may suppose $|J|$ is constant. Fixing a particular $I, \alpha$,
  apply $g_\alpha df_{[n] - I} \odot -$ to \eqref{eq:alt_upstairs.1}.
  If $J \neq I$, then $df_{[n] - I} df_J = 0$ since some $df_i$ appears
  twice. If $J = I$, we have $df_{[n]-I} df_I = \pm df_{[n]} = \pm \Delta \theta_{[n]}$.
  Consequently, we're left with (up to an overall sign)
    \[ 0 = \sum_\beta c_{I, \beta} g_\alpha \Delta \theta_{[n]} \odot g_\beta \Delta
        = \sum_\beta c_{I, \beta} (g_\alpha \Delta \cdot g_\beta \Delta)\psi_{[n]}. \]
  By homogeneity in the $\bm{x}$-variables, we may suppose $\deg g_\alpha
  = \deg g_\beta$, so that $g_\alpha \Delta \cdot g_\beta \Delta = \langle \tilde{g}_\alpha,
  \tilde{g}_\beta \rangle$. By orthogonality of $\{\tilde{g}_\alpha\}$, it follows
  that $c_{I, \alpha} = 0$.
  
  We have just shown
    \[ \Hilb((\bk[V] \otimes \Lambda(V))^{\det}; q, t)
        \geq \prod_{i=1}^n \frac{q^{d_i-1} + t}{1 - q^{d_i}} \]
  where $\geq$ indicates coefficient-wise inequality as bivariate formal power series.
  From \Cref{tab:products}, equality holds, so $\{df_I \odot f^\alpha \Delta\}$ spans
  $(\bk[V] \otimes \Lambda(V))^{\det}$, completing the proof.
\end{proof}

\begin{Remark}
  When $G=S_n$, in the preceding proof we may in fact use
  $\{g_\alpha\} = \{s_\lambda\}$ where
  $s_\lambda$ denotes a Schur polynomial in $n$ variables. More precisely,
  \begin{equation}\label{eq:schur_orthogonal}
    \langle s_\lambda(x_1, \ldots, x_n) \Delta_n, s_\mu(x_1, \ldots, x_n) \Delta_n\rangle
    = (\lambda + \delta_n)! n!\delta_{\lambda, \mu}
  \end{equation}
  where $\Delta_n \coloneqq \prod_{1 \leq i < j \leq n} (x_i - x_j)$ and
  $\delta_n \coloneqq (n-1, n-2, \ldots, 0)$. Indeed, the
  classical bialternant expression \cite[\S7.15]{ec2} gives
  $s_\lambda \Delta_n = \det(x_i^{\lambda_j+n-j}) = \sum_{\sigma \in S_n}
  (-1)^{\ell(\sigma)} \prod_j x_{\sigma(j)}^{\lambda_j+n-j}$. Since
  $\lambda_j+n-j$ is \textit{strictly} decreasing, it follows that when
  $\lambda \neq \mu$, $s_\lambda \Delta_n$ and $s_\mu \Delta_n$
  have no monomials in common. \eqref{eq:schur_orthogonal} now
  follows from \Cref{lem:herm_nondeg}.
\end{Remark}

\section{Harmonic and coinvariant alternants in the real case}\label{sec:alt_harmonics}

Throughout this section,  we assume $\bk \subset \bR$, so $G \leq O(n, \bk)$
consists of orthogonal matrices. Consequently, we may identify the two super-polynomial rings
$\bk[V] \otimes \Lambda(V^*) = \bk[\bm{x}_n, \bm{\theta}_n]$
and $\bk[V] \otimes \Lambda(V) = \bk[\bm{x}_n, \bm{\psi}_n]$ since
$\theta_i \mapsto \psi_i$ is an isomorphism of $G$-modules. In particular,
the two differential operator actions $\cdot$ and $\odot$ now both act on the
same space $\bk[\bm{x}_n, \bm{\theta}_n]$. Since $G$ consists of orthogonal matrices,
$x_1^2 + \cdots + x_n^2$ is $G$-invariant. 

\begin{Definition}
  The \textit{Laplacian} on $\bk[V] \otimes \Lambda(V^*)$ is
    \[ \nabla^2 \coloneqq \sum_{i=1}^n (\partial_i^x)^2. \]
  Thus, $\nabla^2 f = (x_1^2 + \cdots + x_n^2) \cdot f$. 
  By \Cref{thm:unitary_actions}, we then have $\sigma(\nabla^2 f) = \nabla^2 (\sigma(f))$.
  In particular, if $f$ is $G$-invariant, then so is $\nabla^2 f$.
\end{Definition}

The following is an elementary ``polarization identity'' for the Laplacian.

\begin{Lemma}\label{lem:laplacian}
  For all $f, g \in \bk[V]$, we have
    \[ \sum_{i=1}^n \frac{\partial f}{\partial x_i} \frac{\partial g}{\partial x_i}
        = \frac{1}{2}\left(\nabla^2 (fg) - (\nabla^2 f)g - f (\nabla^2 g)\right). \]
  
  \begin{proof}
    By the classical Leibniz rule,
    \begin{align*}
      \frac{\partial^2(fg)}{\partial x_i^2}
        &= \frac{\partial^2 f}{\partial x_i} g + 2\frac{\partial f}{\partial x_i}
             \frac{\partial g}{\partial x_i} + f \frac{\partial^2 g}{\partial x_i}.
    \end{align*}
    Summing over $i=1, \ldots, n$ gives
      \[ \nabla^2(fg) = (\nabla^2 f) g + 2 \sum_{i=1}^n \frac{\partial f}{\partial x_i}
             \frac{\partial g}{\partial x_a} + f (\nabla^2 g). \]
  \end{proof}
\end{Lemma}

\begin{Lemma}
  For $f, g \in \bk[V]$ and $h \in \bk[V] \otimes \Lambda(V^*)$, we have
  \begin{equation}\label{eq:commuting_df}
    df \cdot (dg \odot h) = -dg \odot (df \cdot h) + \frac{1}{2}
    (\nabla^2(fg) - (\nabla^2 f)g - f(\nabla^2 g)) \cdot h.
  \end{equation}
  
  \begin{proof}
  We calculate
    \begin{align*}
      df(\partial^{\bf{x}}, \partial^{\bf{\theta}}) dg(\partial^{\bf{x}}, m^{\bf{\theta}})
        &= \sum_{a=1}^n \frac{\partial f}{\partial x_a}(\partial^{\bm{x}}) \partial_a^\theta
             \sum_{b=1}^n \frac{\partial g}{\partial x_b}(\partial^{\bm{x}}) m_b^\theta \\
        &= \sum_{a, b} \frac{\partial f}{\partial x_a}(\partial^{\bm{x}})
             \frac{\partial g}{\partial x_b}(\partial^{\bm{x}})
             \partial_a^\theta m_b^\theta \\
        &= \sum_{b, a} \frac{\partial g}{\partial x_b}(\partial^{\bm{x}})
             \frac{\partial f}{\partial x_a}(\partial^{\bm{x}})
             (-m_b^\theta \partial_a^\theta + \delta_{a,b}) \\
        &= -\sum_b \frac{\partial g}{\partial x_b}(\partial^{\bm{x}}) m_b^\theta
             \sum_a \frac{\partial f}{\partial x_a}(\partial^{\bm{x}}) \partial_a^\theta
             + \sum_a \frac{\partial g}{\partial x_a}(\partial^{\bm{x}})
                \frac{\partial f}{\partial x_a}(\partial^{\bm{x}}) \\
        &= -dg(\partial^{\bm{x}}, m^{\bm{\theta}}) df(\partial^{\bm{x}},
              \partial^{\bm{\theta}}) + \sum_a \frac{\partial g}{\partial x_a}(\partial^{\bm{x}})
                \frac{\partial f}{\partial x_a}(\partial^{\bm{x}}),
    \end{align*}
    where the third equality follows from \Cref{lem:comm_acomm}.
    The result now follows from \Cref{lem:laplacian}.
  \end{proof}
\end{Lemma}

\begin{Lemma}\label{lem:df_closed}
  Suppose $f \in \bk[V]^G$ is homogeneous with $\deg f \geq 2$.
  Then $df \odot \cH_G \subset \cH_G$.
  
  \begin{proof}
    Suppose $g \in \cH_G$ is harmonic. By \Cref{lem:harmonics_annihilated}, we have
    $f_i \cdot g = df_i \cdot g = 0$ for all $i$, and
    we must show $f_i \cdot (df \odot g) = df_i \cdot (df \odot g) = 0$.
    Since $u \cdot v = u \odot v$ for all $u \in \bk[V]$, we have
      \[ f \cdot (df_i \odot g) = df_i \odot (f \cdot g)
          = df_i \odot 0 = 0. \]
    As for $df_i \cdot (df \odot g)$, by \eqref{eq:commuting_df} we have
    \begin{align*}
      df_i \cdot (df \odot g)
        &= -df \odot (df_i \cdot g) + \frac{1}{2}
             (\nabla^2(f f_i) - \nabla^2(f) f_i - f (\nabla^2 f_i)) \cdot g.
    \end{align*}
    Since $df_i, f_i, f \in \cJ_+^G$, we have
    $df_i \cdot g = f_i \cdot g = f \cdot g = 0$, so each term except
    possibly $\nabla^2(f f_i) \cdot g$ vanishes.
    Since $f, f_i \in \bk[V]$ are $G$-invariant, $\nabla^2(f f_i)$ is also
    $G$-invariant. The result follows trivially if $\nabla^2(f f_i) = 0$, so
    suppose $\nabla^2(f f_i) \neq 0$. Since $\deg f \geq 2$ and
    $\deg f_i \geq 1$, we have $\deg \nabla^2(f f_i) \geq 1$, so
    $\nabla^2(f f_i) \in \cJ_+^G$. Thus indeed $\nabla^2(f f_i) \cdot g = 0$,
    completing the proof.
  \end{proof}
\end{Lemma}

We are now in a position to prove \Cref{thm:alt_harmonics} from the introduction.

\begin{proof}[Proof of \Cref{thm:alt_harmonics}]
  From \Cref{thm:alt_upstairs}, $\{df_I \odot \Delta : I \subset [n]\}$ is
  $\bk$-linearly independent and carries the $\det$-representation. It is well-known that
  $\Delta \in \cH_G$.
  Indeed, $df_i \cdot \Delta = 0$ since $\Delta$ contains no $\theta$'s,
  and $f_i \cdot \Delta$ would be an alternant in $\bk[V]$ of degree lower than
  $\Delta$, so $f_i \cdot \Delta = 0$. By \Cref{lem:df_closed},
  $\{df_I \odot \Delta\} \subset \cH_G$. We must only show this is a spanning set.
  By \Cref{thm:alt_upstairs}, an arbitrary homogeneous element of $\cH_G^{\det} \subset
  \bk[\bm{x}_n, \bm{\theta}_n]^{\det}$ is of the form
    \[ \sum_J df_J \odot g_J \Delta \in \cH_G\]
  for $g_J \in \bk[V]^G$ homogeneous and $|J|$ constant.
  By \Cref{lem:df_closed}, applying $df_{[n] - I}$ for fixed $I$ gives
    \[ \Delta \theta_{[n]} \odot g_I \Delta = (\Delta \cdot g_I\Delta)\theta_{[n]} \in \cH_G. \]
  If $\deg g_I > 0$, then $g_I \in \cJ_+^G$, so applying $g_I \cdot -$ gives
    \[ (g_I\Delta \cdot g_I\Delta)\theta_{[n]} = 0, \]
  forcing $g_I \Delta = 0$, so $g_I = 0$, a contradiction. Hence $g_I$ is constant
  for all $I$, which completes the proof.
\end{proof}

\section{Acknowledgments}

The author would like to thank Vic Reiner, Brendon Rhoades, and Nolan Wallach for
helpful discussions and for sharing their preprints \cite{RSS19,1906.03315,1906.11787};
Mike Zabrocki for sharing the conjecture \cite{1902.08966} which was ultimately
the impetus for the present work; and John Mahacek and Bruce Sagan for helpful
discussions on related work.

\bibliography{refs}{}
\bibliographystyle{acm}

\end{document}